\documentclass[12pt]{article}
\usepackage{graphics, latexsym}
\usepackage{graphicx}
\usepackage{setspace}
\usepackage{amssymb}
\usepackage{amsmath}
\usepackage{epsfig}
\usepackage{CJK}
\usepackage{verbatim}
\usepackage{amsmath, amssymb}
\usepackage{makecell}
\usepackage{multirow}

\begin{document}
\title{ Hypothesis Test of a Truncated Sample Mean for the Extremely Heavy-Tailed Distributions $^*$}
\author{TANG Fuquan $\,\,\,\,\,\,$  HAN Dong$^*$ \\
Department of Statistics, School of Mathematical Sciences,\\
Shanghai Jiao Tong University, Shanghai, 200240, China  }

\maketitle
\begin{abstract}
This article deals with the hypothesis test for the extremely heavy-tailed distributions with infinite mean or variance by using a truncated sample mean.
We obtain three necessary and sufficient conditions under which the asymptotic distribution of the truncated test statistics converges  to normal, neither normal nor stable or converges to $-\infty$ or the combination of stable distributions, respectively. The numerical simulation illustrates an application of the theoretical results above in the hypothesis testing.
\end{abstract}
\renewcommand{\thefootnote}{\fnsymbol{footnote}}
\footnotetext{
${}^*$Supported by National Natural Science Foundation of China (11531001)
\newline ${}^*$ Corresponding author, E-mail: donghan@sjtu.edu.cn }



\section{Introduction}

There are many random systems with the heavy-tailed distribution such as the change in cotton prices, financial market returns, magnitudes of earthquakes and floods, the internet data and so on ([4]).
A random variable with the extremely heavy tails means that its mean or variance is infinite.  Three random variables with  continuous distributions:   Cauchy, Levy and Pareto with  the infinite mean, are often used to model the extremely heavy tails in the stock returns and in the transmission line restoration times (see [8], [9], [16]).

When the mean is finite, there are many methods for  estimating the mean and using the sample mean $\hat{\mu}_n=n^{-1}\sum_{k=1}^n X_k$ to do the hypothesis testing (see [7], [11], [13-17]), where the sequence $\{X_k, k\geq 1\}$ is i.i.d. with the common distribution function $F(x)$. Lugosi and Mendelson [12] gave an overview on mean estimation and regression under heavy-tailed distributions. 

When the variance is infinite, for example, $\textbf{E}|X_k|^{1+\alpha}=\mu_{\alpha}$ for some $\alpha \in (0, 1]$,  Bubeck et al. [2] proposed a truncated empirical mean $\hat{\mu}_T=n^{-1}\sum_{k=1}^nX_kI(|X_k|\leq b_k)$ to estimate the mean $\mu$ and gave the convergence rate of $|\hat{\mu}_T-\mu|$, where $b_k=(bk/\log(1/\delta))^{1/(1+\alpha)}$, $0<\delta<1$,  $b\geq \mu_{\alpha}$ and $I(.)$ is the indicator function. Lee et al. [10] further present a novel robust estimator with the same convergence rate.  Avella-medina et al.[1] considered Huber's M-estimator $\hat{\mu}_H$ of $\mu$ which is defined as the solution to $\sum_{k=1}^n\psi_H(X_k-\mu)=0$, where $\psi_H(x)=\min\{H, \max(-H, x)\}$ is the Huber function. They also obtained the same convergence rate of $|\hat{\mu}_H-\mu|$. 

However, when the mean or variance is infinite, there is not much research work on hypothesis testing. Let $F(x)=1-x^{-\alpha}L(x)$ for $x>0$, where the index $\alpha$ satisfies  $0<\alpha <2$ and $L(x)$ is a positive  slowly vary function. It is well-known (see [3]) that  $n(\hat{\mu}_n-\mu_n)/c_n$ converges in distribution to the stable random variable $S_{\alpha, 1}$ with the index $\alpha \in (0, 2)$ and the shaper parameter 1, where $c_n=\inf\{x: F(x)\geq 1-1/n\}$, $\mu_n=0$ for $0<\alpha <1$ and $\mu_n=\textbf{E}(X_1I(X_1\leq c_n))$ for $1\leq \alpha <2$.  But unfortunately, it is difficult to do the hypothesis testing by using the stable distribution since the stable random variable $S_{\alpha, 1}$ has no analytic expression of density function when $\alpha \neq 1/2, \, 1$ (see [5]).

It is natural to ask a  question: how to construct a test statistic such that it can be easily used to do the hypothesis testing when the mean or  variance is infinite ?
The main purpose of this paper is to solve the problem.

In Section 2,  we present two test statistics with the truncated sample mean. The three necessary and sufficient conditions under which the asymptotic distribution of the truncated sample statistics converges  to normal, neither normal nor stable or converges to $-\infty$ or the combination of stable distributions are given in Section 3.  A relationship between the limiting distribution of the truncated sample statistics and the stable distribution is shown in Section 4.  Section 5 illustrates simulation results about the hypothesis testing for three rejected regions. Section 6 provides some concluding remarks.  Proofs of the theorems are given in the Appendix.

\section{Two test statistics with the truncated sample mean.}
Note that a random variable $X$ can be written as the summation of positive and negative parts, that is, $X=X^+-X^-$. Without loss of generality, we only consider the nonnegative random variables in the following.
Let  $X_k, k\geq 1,$ be a sequence of mutually independent nonnegative random variables with extremely heavy-tailed distribution functions $F_k(x)$.
Similar to the paper [1], we present a  truncated sample mean in the following
\begin{eqnarray}
\hat{\mu}(b_n)=n^{-1}\sum_{k=1}^nX_kI(X_k\leq b_n),
\end{eqnarray}
where $\{b_n\} $ is a positive sequence satisfying $b_n \to \infty$ as $n\to \infty$ which can be called as  the truncated sequence. Note that when $X_kI(X_k\leq b_n)$ is replaced with $\min\{X_k, \, b_n\}$, the resulting  truncated sample mean has the same statistical properties except for a constant factor. 

According to the the sample variance $\hat{B}^2(b_n) $ and the variance $Var(b_n):= \sum_{k=1}^nVar(X_k(b_n))$, we  may present two test statistics
\begin{eqnarray}
T_n:=\frac{n[\hat{\mu}(b_n)-\mu_0(b_n)]}{\hat{B}(b_n)}, \,\,\,\,\,\,\,\,\,\,\,T^o_n:=\frac{n[\hat{\mu}(b_n)-\mu_0(b_n)]}{\sqrt{Var(b_n)}}
\end{eqnarray}
to do the following hypothesis testing:
\begin{eqnarray*}
\text{Original hypothesis} \,\, H_0: \,\, \mu_n=\mu_0(b_n),\,\,\,\,\,\,\,\,\,\,\,             \text{ alternative hypothesis}\,\,  H_1:   \,\, \mu_n\neq \mu_0(b_n),
\end{eqnarray*}
where $\mu_0(b_n)=n^{-1}\sum_{k=1}^n\mu_k(b_n)$, $\mu_k(b_n)=\textbf{E}(X_k(b_n))$, $X_k(b_n)=X_kI(X_k\leq b_n)$ for $k\geq 1$ and $\hat{B}(b_n)=[\sum_{k=1}^n(X_k(b_n)-\hat{\mu}(b_n))^2]^{1/2}$.

This leads to another problem: how to choose the truncated sequence  $\{b_n\}$ such that the distributions of the test statistics $T_n$ and $T^o_n$ can converge to some limiting distribution? We are particularly concerned that how to choose the truncated sequence  $\{b_n\}$ such that $T_n$ and $T^o_n$ can converge to normal distribution.

To this end, we first present two conditions. Let  $X_k, 1\leq k\leq n,$ be mutually independent non-negative random variables with the extremely heavy-tailed distribution functions $F_k(x)$ satisfying
\begin{eqnarray*}
(I)\,\,\,\,\,\,\,\,\,\,\,\,\,\,\,\,\,\,\,\,\,\,\,\,\,\,\,\,\,\,\,\,\,\,\,\,\,\,\,\,\,\,\,\,\max_{1\leq k\leq n}\{1-F_k(b_n)\}\to 0  \,\,\,\,\,\,\,  ( n\to \infty )
\end{eqnarray*}
and there is  a series of positive numbers $d_k(r)=\alpha_k/(r-\alpha_k)$ satisfying $0<\alpha_k<1, \alpha_k \to \alpha, 0<\alpha<1$ such that
\begin{eqnarray*}
(II)\,\,\,\,\,\,\,\,\,\,\,\,\,\,\,\,\,\,\,\,\,\, \textbf{E}(X^{r}_k(b_n))=(1+o(1))d_k(r)b^{r}_n[1-F_k(b_n)] \to \infty  \,\,\,\,\,\,  ( n\to \infty )
\end{eqnarray*}
for $r\geq 1$ and $1\leq k\leq n$.

There are many probability density functions, for example, the following five density functions, satisfying the two conditions (I) and (II).
\begin{eqnarray*}
f_k(x)&=&\frac{c_{k1} cos(x^{-\alpha_k})}{x^{1+\alpha_k}(sin(x^{-\alpha_k})+1)},\,\,\,\,\, g_k(x)=c_{k2}x^{-(\alpha_k+1)}cos(x^{-1}),\\
h_k(x)&=&c_{k3}x^{-\alpha_k}sin(1/x), \,\,\,\,\, p_k(x)=c_{k4}sin(x^{-(\alpha_k+1)}), \,\,\,\,\, q_k(x)=c_{k5} x^{-(\alpha_k+1)}
\end{eqnarray*}
for $x\geq 1$ and $k\geq 1$,  where $c_{kl}, 1\leq l\leq 5,$ are  the normalized positive numbers, $0<\alpha_k<1$ and $\alpha_k \to \alpha, 0<\alpha<1$. In fact, for $\{X_k\}$ with density functions $\{f_k(x)\}$ we have
\begin{eqnarray*}
1-F_k(b_n)=\frac{c_{k1}}{\alpha_k}\log (sin(b^{-\alpha_k}_n)+1)=(1+o(1))\frac{c_{k1}}{\alpha_k}b^{-\alpha_k}_n\, \to \, 0
\end{eqnarray*}
uniformly for $k\geq 1$ as $n\to \infty$ and
\begin{eqnarray*}
\textbf{E}(X^r_k(b_n))=(1+o(1))\int_{M}^{b_n}x^rf_k(x)dx=(1+o(1))\frac{c_{k1}b_n^{r-\alpha_k}}{r-\alpha_k}=(1+o(1))\frac{\alpha_k}{r-\alpha_k} b_n^r[1-F_k(b_n)]
\end{eqnarray*}
for $r\geq 1$ and large $M$, where $b_n\gg M$. That is, the conditions (I) and (II) hold for $\{f_k(x)\}$. Similarly,  we can check that the other four kinds of probability densities   $\{g_k(x)\}, \{h_k(x)\}, \{p_k(x)\}$ and $\{q_k(x)\}$ also satisfy the two conditions.

Next we consider the case that $1\leq \alpha_k<2$. For this case, we may have $\textbf{E}(X_k)<\infty$  and $\textbf{E}(X^2_k)=\infty$. In addition to the condition (I),  the condition (II) for this case will be replaced by the following condition (II)$^{\prime}$
\begin{eqnarray*}
(II)' \,\,\,\,\,\,\,\,\,\,\,\,\,\,\,\,\,\,\,\,\, \textbf{E}(X^{r}_k(b_n))=(1+o(1))d_k(r)b^{r}_n[1-F_k(b_n)]\to \infty  \,\,\,\,\,\,  ( n\to \infty )
\end{eqnarray*}
for $r\geq  2$ and $1\leq k\leq n$, where  $d_k(r)=\alpha_k/(r-\alpha_k)$ for $r\geq 2$ satisfying $1\leq \alpha_k<2, \alpha_k \to \alpha, 1\leq \alpha<2$.

For the  five kinds of  probability densities   $f_k(x), g_k(x), h_k(x), p_k(x)$ and $q_k(x)$ above, with $1\leq \alpha_k<2, \alpha_k \to \alpha, 1\leq \alpha<2$,  we can similarly check that they all satisfy the conditions (I) and (II)$^{\prime}$.

\textbf{Remark 1.} Note that the probability density $u_{\alpha_k}(x)$ of the stable random variable $S_{\alpha_k, 1}$ can be written as (see [5]) $u_{\alpha_k}(x)=c_k(1+o(1))x^{\alpha_k+1}$ for large $x$, where $0<\alpha_k<2, k\geq 1$.  We can also check that  $u_{\alpha_k}(x)$ with $0<\alpha_k<1, \alpha_k \to \alpha, 0<\alpha<1,$ and  $u_{\alpha_k}(x)$ with $1\leq \alpha_k<2, \alpha_k \to \alpha, 1\leq \alpha<2$, satisfy the conditions (I), (II) or  (II)$^{\prime}$, respectively.

The following lemma shows that the regularly varying distribution functions (see [2]) also satisfies the conditions (I), (II) or (II)$^{\prime}$.

Let the regularly varying distribution functions $F_k(x), k\geq 1,$ satisfy
\begin{eqnarray}
1-F_k(x)=x^{-\alpha_k}L(x)\,\,\,\,\,\,\, 0< \alpha_k<1 \,\,\text{ or}\,\,\,1\leq  \alpha_k<2
\end{eqnarray}
for $x\geq a_k>0$, where $\{a_k\}$ is bounded, $\alpha_k \to \alpha$, ($0<\alpha<1$ or $1\leq \alpha<2$) and $L(x)$ is a positive monotone continuous slowly vary function. Here, the slowly vary function means that
\begin{eqnarray}
L(tx)/L(t) \to 1
\end{eqnarray}
for any positive number $x$ as $t\to \infty$. For example, $L(x)=\theta(\log(1+x))^{\tau}$, $L(x)=\log\log x, \, x>e$ and their reciprocals, are all the slowly vary functions, where $\tau, \theta$ are two positive constants.

\textbf{Lemma 1.} \textit{ The non-negative random variables $X_k, k\geq 1,$ with the regularly varying distribution functions $F_k(x), k\geq 1,$ satisfying (3) and (4) also satisfy the conditions (I), (II) or  (II)$^{\,\prime}$.}

\section{ Asymptotic distributions of the two test statistics}

\subsection{ Asymptotic distributions of $T_n$ }
Let
\begin{eqnarray}
\xi_n:&=&\frac{\sum_{k=1}^n(X_k(b_n)-\mu_k(b_n))}{b_n}, \,\,\,\,\,\,\,\,\,\,\,\eta_n:=\frac{\sum_{k=1}^n(X_k(b_n)-\mu_0(b_n))^2}{b^2_n},\\
h_n:&=&\sum_{k=1}^n(1-F_k(b_n)), \,\,\,\,\,\,\,\,\,\,\,\,\,\,\,\,\,\,\,\,\,\,   h_n(m):=\sum_{k=1}^nd_k(m)(1-F_k(b_n))
\end{eqnarray}
for $m\geq 1$ and $n\geq 1$. Let $N(0, 1)$ and $\Rightarrow$  denote the standard normal distribution and the convergence in distribution, respectively.

The following  theorem gives the three sufficient conditions, $h_n \to \infty$, $h_{n}\to \, h$ and $h_n \to 0$ for the asymptotic convergence  of $T_n$.

\textbf{Theorem 1.} \textit{  Assume that the original hypothesis $H_0$ holds and $X_k, k\geq 1,$ are mutually independent non-negative random variables with the extremely heavy-tailed distribution functions $\{F_k(x)\}$ satisfying the conditions (I) and (II) or (II)$^{\,\prime}$. \\
(i) If $h_n \to \infty$, then $T_n \Rightarrow N(0, 1)$;\\
(ii) If  $h_{n}\to \, h$  for some positive number $h$, then, there are two random variables $\xi, \eta$
such that $ (\xi_{n},\,\eta_{n}) \Rightarrow (\xi,\, \eta)$ and $ T_{n} \Rightarrow \xi/\sqrt{\eta}$, where both $\xi$ and $\eta$ have the following characteristic  functions
\begin{eqnarray}
C_{\xi}(t)&=&\exp\{h\alpha |t|^{\alpha}\Big(i\, sgn(t)\int_{0}^{|t|}\frac{sinx-x}{x^{1+\alpha}}dx-\int_{0}^{|t|}\frac{1-cosx}{x^{1+\alpha}}dx\Big)\}\\
C_{\eta}(t)&=&\exp\{h\frac{\alpha}{2} |t|^{\frac{\alpha}{2}}\Big(i\, sgn(t)\int_{0}^{|t|}\frac{sinx}{x^{1+\frac{\alpha}{2}}}dx-\int_{0}^{|t|}\frac{1-cosx}{x^{1+\frac{\alpha}{2}}}dx\Big)\}
\end{eqnarray}
respectively for $0<\alpha <2$. \\
(iii) If $h_n \to 0$ and $0<\alpha\leq 1$, then  $T_n \to -\infty$ in probability;\\
(iv) If $h_n \to 0$ and $1< \alpha<2$, then,  $(\xi_{n}h_n^{-1/\alpha},\,\,\eta_{n}h_n^{-2/\alpha}) \Rightarrow (S_{\alpha, -1},\,\, S_{\alpha/2, 1})$ and $T_{n} \Rightarrow S_{\alpha, -1}/\sqrt{S_{\alpha/2, 1}}$, where  both $S_{\alpha, -1}$ and $S_{\alpha/2, 1}$ are stable random variables with the shaper parameters $-1$ and $1$, respectively.}

\textbf{Remark 2.}  It can be seen from the proof of (i) of Theorem 1 that the result in (i) still holds if we use the following weak conditions:  there is a positive sequence $\{d_k(r)\}$ with
\begin{eqnarray*}
0<\min_{k\geq 1}\{d_k(r)\}\leq d_k(r) \leq \max_{k\geq 1}\{d_k(r)\}<\infty
\end{eqnarray*}
satisfies the conditions (II) or (II)$^{\,\prime}$ for $r=1, 2, 3, 4$ or $r=2, 3, 4$, respectively.  And the conditions  that $d_k(r)=\alpha_k/(r-\alpha_k)\to \alpha/(r-\alpha)$ ($0<\alpha<1$ or $1\leq \alpha<2$) satisfying the conditions (II) or (II)$^{\,\prime}$, can let us to obtain the display expressions of the four characteristic  functions of random variables $\xi,\, \eta$, $S_{\alpha, -1}$ and $S_{\alpha/2, 1}$.

\textbf{Remark 3.} The distribution functions of $\xi$ and $S_{\alpha, -1}$ are not symmetric, so do  $\xi/\sqrt{\eta}$ and $S_{\alpha, -1}/\sqrt{S_{\alpha/2, 1}}$.  In other words, the distributions of  both  $\xi/\sqrt{\eta}$ and $S_{\alpha, -1}/\sqrt{S_{\alpha/2, 1}}$ are not the standard normal, $N(0, 1)$. Moreover, $\xi, \eta$ have continuous probability density functions since the absolute value functions of the characteristic  functions of  $\xi, \eta$ are integrable.

It is clear that $h_n \to \infty$ or $h_n \to 0$ mean that the truncated sequence $\{b_n\}$ approaches  to infinity slowly or goes  to infinity quickly, respectively.

As an application of (i) of Theorem 1, we can get the confidence interval  for $\mu_0(b_n)$ in the following
\begin{eqnarray}
\textbf{P}\Big(   \mu_0(b_n)\in [ \hat{\mu}(b_n)-\frac{x\hat{B}(b_n)}{n},\,\,\,\,\frac{x\hat{B}(b_n)}{n}+\hat{\mu}(b_n)]\Big)\approx 2\Phi(x)-1\,\,\,\,\, (\, x\,>\,0 \,)
\end{eqnarray}
for large $n$ with the probability $2\Phi(x)-1$, where $\Phi(x)$ is the standard normal distribution.  It can be seen that if $b'_n > b_n$ with $\sum_{k=1}^n[1-F_k(b'_n)] \to \infty$ and $\sum_{k=1}^n[1-F_k(b_n)] \to \infty$, then the mean $\mu_0(b'_n)$ is closer the real mean than the mean $\mu_0(b_n)$.

By Theorem 1 we can further get the following theorem which gives the three  necessary and sufficient conditions of asymptotic convergence of $T_n$.

\textbf{Theorem 2.} \textit{  Let the conditions of the theorem 1 hold. \\
(i) $T_n \Rightarrow N(0, 1)$  if and only if $h_n \to \infty$;\\
(ii) $T_{n}\Rightarrow \xi/\sqrt{\eta}$ if and only if $h_{n} \to \, h \,\,\, (\,0<h<\infty\,)$. \\
(iii) Let $0<\alpha\leq 1$.  $T_{n}\to -\infty$ in probability if and only if   $h_{n} \to 0$.\\
(iv) Let $1<\alpha<2$.  $T_{n}\Rightarrow S_{\alpha, -1}/\sqrt{S_{\alpha/2, 1}}$ if and only if  $h_{n} \to 0$.}

Next, we discuss on the set of  the truncated sequence  $\{b_n\}$.  Let $D_n(x)=\sum_{k=1}^n[1-F_k(x)]$ and define three regions or sets $R_s, R_b$ and $R_c$ of the truncated sequence  $\{b_n\}$   in the following
\begin{eqnarray*}
R_s=\{\{b_n\}:\,D_n(b_n) \to \infty\}, \,\, R_c=\{\{b_n\}:\,h_1\leq D_n(b_n)\leq h_2\} \,\text{ and }\, R_b=\{\{b_n\}:\,D_n(b_n) \to 0\},
\end{eqnarray*}
where $h_1, h_2$ are two positive numbers. Let the truncated sequence $\{c_n\} \in R_c$. We may call $\{c_n\}$ and $R_c$  as a critical truncated sequence and the critical truncated region respectively, since when the truncated sequence $\{b_n\}\in R_s$, that is, $\{b_n\}$ is smaller  than the critical truncated sequence $\{c_n\}$,  then $T_n \Rightarrow N(0, 1)$, when the truncated sequence $\{b_n\}\in R_b$, that is, $\{b_n\}$ is bigger than the critical truncated sequence $\{c_n\}$,  then   $T_{n}$ converges  to $-\infty$ in probability or $T_{n}\Rightarrow S_{\alpha, -1}/\sqrt{S_{\alpha/2, 1}}$, and for the critical truncated sequence, there is a subseries $\{n_l\}$ such that $T_{n_l}\Rightarrow \xi/\sqrt{\eta}$.

Now we consider an application of Theorem 1 to the distribution functions $\{F_k(x)\}$ satisfying (3) and (4). Let $l(n)$ is an increasing positive sequence satisfying $l(n)/n \to 0$ as $n\to \infty$. Assume that
\begin{eqnarray}
(\alpha_k-\alpha)\log n \to 0
\end{eqnarray}
for $k\geq l(n)$ as $n\to \infty$. Take the truncated sequence $\{c_n\} $ satisfying $c_n=(nL(c_n)/h)^{1/\alpha}$ for $0<h<\infty$, it follows that
\begin{eqnarray*}
\sum_{k=1}^n[1-F_k(c_n)] &=&\frac{h}{n}\sum_{k=1}^{n}\Big(\frac{h}{nL(c_n)}\Big)^{\frac{\alpha_k-\alpha}{\alpha}}\\
&=&(1+o(1))\frac{n-l(n)}{n}\frac{h}{n-l(n)} \sum_{k=l(n)}^{n}e^{\frac{(\alpha_k-\alpha)}{\alpha}(\log h -\log n-\log L(c_n))}\\
&=&(1+o(1))h
\end{eqnarray*}
for large $n$. This means that $\{c_n\} \in R_c $ is the critical truncated sequence. It is clear that  $\{b_n\} \in R_s $ for $b_n\leq c_n/\widetilde{L}(n)$ and $\{b_n\} \in R_b $ for $b_n\geq c_n\widetilde{L}(n)$, where $\widetilde{L}(n)$ is a positive monotone continuous slowly vary function satisfying $\widetilde{L}(n)\to \infty$ as $n\to \infty$.

By Lemma 1 and Theorem 1 we have the following corollary.

\textbf{Corollary 1.} \textit{  Let  $X_k, 1\leq k\leq n,$ be mutually independent non-negative random variables with the distribution functions $\{F_k(x)\}$ satisfying the conditions (3), (4) and (10). \\
(i) If the truncated sequence $\{b_n\} $ satisfies $b_n\leq c_n/\widetilde{L}(n)$, then $T_n\Rightarrow N(0, 1)$.\\
(ii) For the critical truncated sequence $\{c_n\} $, $T_{n}\Rightarrow \xi/\sqrt{\eta}$.\\
(ii) Let  $0<\alpha\leq 1$. If $\{b_n\} $ satisfies $b_n\geq c_n\widetilde{L}(n)$, then  $T_{n}$ converges to $-\infty$ in probability.\\
(iii) Let  $1<\alpha<2$. If $\{b_n\} $ satisfies $b_n\geq c_n\widetilde{L}(n)$, then $T_{n}\Rightarrow S_{\alpha, -1}/\sqrt{S_{\alpha/2, 1}}$.}

\subsection{ Asymptotic distributions of $T^o_n$ }

The following theorem  gives three  necessary and sufficient conditions of asymptotic convergence of $T^0_n$.

\textbf{Theorem 3.} \textit{  Under the conditions of Theorem 1, we have \\
(i) $T^o_n \Rightarrow N(0, 1)$  if and only if $h_n \to \infty$;\\
(ii) $T^o_n \Rightarrow T_{\alpha, h}$ if and only if $h_n \to h, 0<h<\infty$,
where  the characteristic function of random variable $T_{\alpha, h}$ has the following form
\begin{eqnarray}
C_{\alpha, h}(t)=\exp\{h\alpha \sigma^{-\alpha}|t|^{\alpha}\Big(i\, sgn(t)\int_{0}^{|t|/\sqrt{h}\sigma}\frac{sinx-x}{x^{1+\alpha}}dx-\int_{0}^{|t|/\sqrt{h}\sigma}\frac{1-cosx}{x^{1+\alpha}}dx\Big)\}
\end{eqnarray}
for $0<\alpha <2$ and $\sigma=\sqrt{\alpha}/\sqrt{2-\alpha}$;\\
(iii)  $T^o_n $ converges  to $0$ in probability  if and only if $h_n \to 0$.}

It can be seen that the distribution of $T_{\alpha, h}$ is neither normal nor stable. So we may call the distribution of $T_{\alpha, h}$  NNS-distribution. 

\textbf{Remark 4.} Similar to (iv) of Theorem 1 we can check that $T_{\alpha, h}/h^{1/\alpha} \Rightarrow S_{\alpha, -1}$ as $h\to 0$ for $1<\alpha <2$. In other words, though the random variable $T_{\alpha, h}$  is not the stable random variable, its some limiting is.

\textbf{Remark 5.} From (iii)-(iv) of Theorem 1 and (iii) of Theorem 3 it follows that the asymptotic statistical properties of $T_n$ and $T^o_n$ are completely different. The reason is that $T^o_n/T_n \to 0$ in probability when $h_n\to 0$. In fact, for small $\varepsilon>0$ and $0<\alpha <1$, by the condition (II) we have
\begin{eqnarray*}
\textbf{P}\Big(\frac{T^o_n}{T_n}\geq \varepsilon\Big)&\leq &\frac{\textbf{E}(\hat{B}(b_n))}{\varepsilon \sqrt{Var(b_n)}}\leq \frac{\sum_{k=1}^n\textbf{E}|X_k(b_n)-\hat{\mu}(b_n)|}{\varepsilon \sqrt{Var(b_n)}}\\
&\leq &\frac{2\sum_{k=1}^n\textbf{E}(X_k(b_n))/b_n}{\varepsilon \sqrt{Var(b_n)/b^2_n}}\leq O(\sqrt{h_n})\to 0
\end{eqnarray*}
as $h_n\to 0$, where  $Var(b_n)=(1+o(1))\alpha_k(2-\alpha_k)^{-1}b^2_n(1-F_k(b_n))$. Note that $\hat{B}^2(b_n)/b^2_n=\eta_n-\xi^2_n/n$. For $1<\alpha <2$, it follows from (iv) of Theorem 1 that
\begin{eqnarray*}
\textbf{P}\Big(\frac{T^o_n}{T_n}\geq \varepsilon\Big)&= &\textbf{P}\Big(\frac{\eta_n-\xi^2_n/n}{Var(b_n)/b^2_n}\geq \varepsilon^2\Big)\\
&=&\textbf{P}\Big(h^{2/\alpha}_n\frac{[\eta_nh^{-2/\alpha}_n-(\xi_nh^{-1/\alpha}_n)^2/n]}{O(h_n)}\geq \varepsilon^2\Big)\\
&=&\textbf{P}\Big((1+o(1))h^{\frac{2-\alpha}{\alpha}}_n[\eta_nh^{-2/\alpha}_n-(\xi_nh^{-1/\alpha}_n)^2/n]\geq \varepsilon^2\Big) \to 0
\end{eqnarray*}
as $h_n\to 0$ since $(\xi_{n}h_n^{-1/\alpha},\,\,\eta_{n}h_n^{-2/\alpha}) \Rightarrow (S_{\alpha, -1},\,S_{\alpha/2, 1})$.

Next we consider  the general  truncated sequence  $\{b_{nk}\}$ satisfying $b_{nk}\to \infty$ for every $k\geq 1$ as $n\to \infty$.  Similar to $T^o_n=T^o_n(b_n)$  we can define the test statistics $T^o_n({b_{nk}})$ for $n\geq 1$ in the following
\begin{eqnarray}
T^o_n({b_{nk}})=\frac{n}{\sqrt{Var(b_{nk})}}\Big(\hat{\mu}(b_{nk})-\mu_0(b_{nk})\Big),
\end{eqnarray}
where $\hat{\mu}(b_{nk})=n^{-1}\sum_{k=1}^nX_k(b_{nk}), \mu_0(b_{nk})=n^{-1}\sum_{k=1}^n\mu_k(b_{nk})$, $\mu_k(b_{nk})=\textbf{E}(X_k(b_{nk}))$, $X_k(b_{nk})=X_kI(X_k\leq b_{nk})$ for $k\geq 1$ and $Var(b_{nk}):= \sum_{k=1}^nVar(X_k(b_{nk}))$.

Let $F_k(x)$ be continuous with $1-F_k(x)=x^{-\alpha_k}L(x)$ and $c_{nk}=\inf\{x: F_k(x)\geq 1-h/n\}$ for $n>h, 1\leq k\leq n.$ Hence, $1-F_k(c_{nk})=h/n$ and $c_{nk}=n^{1/\alpha_k}\tilde{L}_k(n)$, where $\tilde{L}_k(n)=[L(c_{nk})/h]^{1/\alpha_k}$. Similar to Theorem 3, we can get the following theorem 4 which shows that $\{c_{nk}\}$ is a critical truncated sequence.

\textbf{Theorem 4.} \textit{  Let $X_k, 1\leq k\leq n,$ be mutually independent non-negative random variables with the continuous distribution functions $F_k(x), k\geq 1,$ satisfying  (3), (4) and (10).  Let the variance $Var(b_{nk})$ be known and $L_0(x)$ be a positive slowly vary function satisfying $L_0(x)\to \infty$ as $n\to \infty$. Under the original hypothesis $H_0: \mu_n=\mu_0(b_{nk})$, we have that \\
(i) if $b_{nk}\leq c_{nk}L^{-1}_0(c_{nk})$, then,  $T^o_n({b_{nk}})\Rightarrow N(0, 1)$;\\
(ii) if $b_{nk}\geq c_{nk}L_0(c_{nk})$,  then,  $T^o_n({b_{nk}})$ converges to $0$ in probability;\\
(iii) if $b_{nk}= c_{nk}$, then,  $T^o_n(b_{nk})\Rightarrow T_{\alpha, h}$.}

By (iii) of Theorem 4 we see that for $\{b_{nk}\}\in R_c(M):=\{\{b_{nk}\}:\, M^{-1}c_{nk}\leq b_{nk}\leq Mc_{nk}\}$ for any large number $M$, there is a subsequence $\{n_l\}$ such that  $T^o_{n_l}(b_{n_{l}k})\Rightarrow T_{\alpha, h'}$ as $l\to \infty$, where $0<h'<\infty$.

\section{The relationship between $T_{\alpha, 1}$ and $S_{\alpha, 1}$}

By the remark 4, we know that there is a relationship between $T_{\alpha, h}/h^{1/\alpha}$ and $S_{\alpha, -1}$ for $1<\alpha<2$. In this section, we will show another relationship between $T_{\alpha, 1}$ and  $S_{\alpha, 1}$.

Let $X_k, 1\leq k\leq n,$ be i.i.d. positive random variables with the heavy-tailed continuous distribution function $F(x)=1-x^{-\alpha}L(x)$ satisfying $0<\alpha <1$. We have $c_n=\inf\{x:\, F(x)\geq 1-1/n\}=n^{1/\alpha}\tilde{L}(n)$ or $c^{-\alpha}_nL(c_n)=1/n$ and therefore, $D_n(c_n)=\sum_{k=1}^n[1-F(c_n)]=1$,  where $\tilde{L}(n)=[L(c_{n})]^{1/\alpha}$. It is well-known that (see [3])
\begin{eqnarray*}
S_{n, \alpha, 1}:=c^{-1}_n\sum_{k=1}^n(X_k-\mu_n) \Rightarrow S_{\alpha, 1}
\end{eqnarray*}
as $n\to \infty$, where $\mu_n=0$ for $0<\alpha<1$ and $\mu_n=\textbf{E}(X_1(c_n))$ for $1\leq \alpha <2$.

By the critical  truncated sequence $\{c_n\}$, we may decompose $S_{n, \alpha, 1}$ into two parts, $S_{n, \alpha, 1}=U_n+V_n$,  where  $U_n$ and $V_n$ are defined in the following
\begin{eqnarray}
U_n=c^{-1}_n\sum_{k=1}^n[X_kI(X_k\leq c_n)-\mu_n],\,\,\,\,\,\,\,\,\,\,\, V_n=c^{-1}_n\sum_{k=1}^nX_kI(X_k>c_n)
\end{eqnarray}
for $n\geq 1$. In fact, both $U_n$ and $V_n$ can be considered as the sum of the lower parts of $S_{n, \alpha, 1}$ and the sum of the upper parts of $S_{n, \alpha, 1}$, respectively.

The following theorem shows that the limiting distribution of  $U_n$ is the distribution of $T_{\alpha, 1}$ and the limiting distribution of  $V_n$ is the distribution of linear  combination of $T_{\alpha, 1}$ and $S_{\alpha, 1}$.

\textbf{Theorem 5.} \textit{ Suppose $X_k, k\geq 1,$ are i.i.d. positive random variables  with the heavy-tailed continuous distribution function $F(x)=1-x^{-\alpha}L(x)$ satisfying $0<\alpha <2$. Then\\
(i)
\begin{eqnarray*}
U_n \Rightarrow \, \sqrt{\frac{\alpha}{2-\alpha}}\,T_{\alpha, 1}+\frac{\alpha}{1-\alpha},\,\,\,\,\,\,\,\,\,\,\,\,\,\,\,V_n \, \Rightarrow \, S_{\alpha, 1}-\sqrt{\frac{\alpha}{2-\alpha}}\,T_{\alpha, 1}+\frac{\alpha}{1-\alpha}
\end{eqnarray*}
for $0<\alpha<1$;\\
(ii)
\begin{eqnarray*}
U_n \Rightarrow \, \sqrt{\frac{\alpha}{2-\alpha}}\,T_{\alpha, 1},\,\,\,\,\,\,\,\,\,\,\,\,\,\,\,\,\,\,\,\,V_n \, \Rightarrow \, S_{\alpha, 1}-\sqrt{\frac{\alpha}{2-\alpha}}\,T_{\alpha, 1}
\end{eqnarray*}
for $1\leq \alpha<2$.}

\section{ Numerical simulations}

In this section, as an application of (i) of the theorems 1 and 3, the hypothesis testing will be studied by the numerical simulations based on the asymptotic distribution of $T_n$ and $S_{n, \alpha, 1}$.

Let  $X_k, k\geq 1,$ be i.i.d.  with extremely heavy-tailed distribution function $F(x)=1-1/x^{\alpha_0}$ for $x\geq 1$, where $0<\alpha_0<2$. It is clear that $c_n=n^{1/\alpha_0}$. Let the  truncated sequence  $\{b_n\}$ satisfy $b_n\leq c_n/\log (1+n)=n^{1/\alpha_0}/\log (1+n)$. By (i) of the theorems 1 and 3  we know that  $T_n\Longrightarrow N(0, 1)$ and $T^o_n\Longrightarrow N(0, 1)$. Since
\begin{eqnarray*}
\mu_0(b_n)=\frac{\alpha_0}{1-\alpha_0}\Big(b_n^{1-\alpha_0}-1\Big),\,\,\,\,\,\,\,\,\,\,\,\, Var(X_1(b_n))=\frac{\alpha_0}{2-\alpha_0}\Big(b_n^{2-\alpha_0}-1\Big)-\mu^2_0(b_n),
\end{eqnarray*}
it follows that the original hypothesis and the alternative hypothesis  for the truncated mean $\mu_0(b_n)$ can be equivalently written as $H_0: \alpha =\alpha_0$ and $H_1: \alpha\neq \alpha_0$, respectively.

We know that  $S_{\alpha, 1}$ has no analytic expression of density function for $0<\alpha <1$ except $\alpha=1/2$. Let $\alpha=1/2$, we have $c_n=n^2$ and the stable random variable $S_{1/2, 1}$ has the density function (see [5] P.165)
\begin{eqnarray*}
p(x, 1/2)=(2\pi x^3)^{-\frac{1}{2}}e^{-\frac{1}{2x}}  \,\, \text{ for } x >0.
\end{eqnarray*}
In this case, by using the sample mean $\hat{\mu}_n=n^{-1}\sum_{k=1}^nX_k$ we can get the rejected region $\widetilde{R}_{\beta}(n)$ of the original hypothesis $H_0: \, \alpha=\alpha_0=1/2$ with the small probability $\beta=0.05$  in the following
\begin{eqnarray*}
\widetilde{R}_{\beta}(n)=\Big(ny,\,\,\,+\infty\Big),
\end{eqnarray*}
where $y=225$ satisfies
\begin{eqnarray*}
&&\textbf{P}\Big(\frac{\hat{\mu}_n}{n}>y\,|\,H_0\Big)=\int_{y}^{+\infty} p(x, 1/2)dx =0.05
\end{eqnarray*}
for large $n$, since $S_{n, 1/2, 1}=\hat{\mu}_n/n=\sum_{k=1}^nX_k/n^2 \, \Longrightarrow S_{1/2, 1}$. That is, $\textbf{P}(\hat{\mu}_n\in \widetilde{R}_{\beta}(n)\,|\,H_0)\approx 0.05$ for large $n$.
 
From $2\Phi(x)-1 =1-\beta$ it follows that $\Phi(x_{\beta})=1-\beta/2$ and therefore, $x_{\beta}=1.96$ for $\beta=0.05$.  Hence, by (i) of Theorem 3 we can get  the rejected region of the original hypothesis $\alpha=\alpha_0$ with the small probability $\beta=0.05$ and the known variance $Var(X_1(b_n))$ for the truncated sample mean $\hat{\mu}(b_n)$ in the following
\begin{eqnarray*}
R^o_{\beta}(b_n)=\Big(-\infty,\,\,\,\mu_0(b_n)-\frac{1.96 \sqrt{ Var(b_n)}}{n}\Big)\cup \Big( \mu_0(b_n)+\frac{1.96  \sqrt{Var(b_n)}}{n},\,\,\,+\infty\Big),
\end{eqnarray*}
that is, $\textbf{P}(\hat{\mu}(b_n)\in R_{\beta}(b_n)|\,H_0)\approx 0.05$ for large $n$, where $Var(b_n)=nVar(X_1(b_n))$. Similarly, by (i) of theorem 1 we can obtain the rejected region of the original hypothesis $\alpha=\alpha_0$ with the small probability $\beta=0.05$ and the unknown variance for the truncated sample mean $\hat{\mu}(b_n)$ in the following
\begin{eqnarray*}
R_{\beta}(b_n)=\Big(-\infty,\,\,\,\mu_0(b_n)-\frac{1.96 \hat{B}(b_n)}{n}\Big)\cup \Big( \mu_0(b_n)+\frac{1.96 \hat{B}(b_n)}{n},\,\,\,+\infty\Big),
\end{eqnarray*}
where $\hat{B}(b_n)=[\sum_{k=1}^n(X_k(b_n)-\hat{\mu}(b_n))^2]^{1/2}$.

Next we will give the numerical simulations in the following Tables 1 and 2 for the original hypothesis $H_0: \, \alpha=\alpha_0=1/2$ by taking different total numbers of samples,  $n=10^3, 10^4, 10^5$ and the different truncated sequences: $b_n$ is equal to $\log(n)$, $n^{0.5}$,  $n$, $n^{4/3}/10\log(n)$, $n^{1.5}$, $n^{1.8}$ and $n^{2}/\log (1+n)$. Let $N=10^4, N^5$ be the numbers of repetitions and let   $\widetilde{r}_N(n),  r^o_{N}(b_n)$ and $r_{N}(b_n)$ denote the percentages of times that the sample mean or the truncated sample mean falls in the three rejected regions, $\widetilde{R}_{\beta}(n),  R^o_{\beta}(b_n)$ and $R_{\beta}(b_n)$, respectively.

The percentages $r^o_{N}(b_n)$, $r_{N}(b_n)$ and $\widetilde{r}_N(n)$ are listed in penultimate columns and last columns in Tables 1 and 2, respectively. 

\begin{center}
\tabcolsep 0.0000000001in
\setlength{\tabcolsep}{6mm}{
\vskip 0.3cm
{{\bf Table 1:} Simulation  values of  $r^o_{N}(b_n)$ and $\widetilde{r}_N(n)$ for $\alpha_0=1/2$, $\beta=0.05$  and $y=255.$   
 \\}
\vskip 0.6cm
\centering
\begin{tabular}{|c|c|c|c|c|c|}\hline

$b_n$                                           &$n$      & $N$     & $r^o_{N}(b_n)$    &$\widetilde{r}_N(n)$    \cr\cline{1-5}
\multirow{3}{*}{$\log(n)$}                      & $10^3$  & $10^4$  & $0.0492$                              & $0.0442$                          \\\cline{2-5}
		                                        & $10^4$  & $10^4$  & $0.0505$                              & $0.0507$                          \\\cline{2-5}
                                                & $10^5$  & $10^5$  & $0.0502$                             & $0.04983$                         \\\hline

\multirow{3}{*}{$n^{0.5}$}                      & $10^3$  & $10^4$  & $0.0488$                              & $0.0405$                          \\\cline{2-5}
		                                        & $10^4$  & $10^4$  & $0.0498$                              & $0.0463$                          \\\cline{2-5}
                                                & $10^5$  & $10^5$  & $0.0498$                             & $0.0501$                         \\\hline

\multirow{3}{*}{$n $}                           & $10^3$  & $10^4$  & $0.0471$                              & $0.0496$                          \\\cline{2-5}
		                                        & $10^4$  & $10^4$  & $0.0505$                              & $0.0507$                          \\\cline{2-5}
                                                & $10^5$  & $10^5$  & $0.0517$                             & $0.0499$                         \\\hline

\multirow{3}{*}{$\frac{n^{4/3}}{10\log(n)}$}   & $10^3$  & $10^4$  & $0.0534$                              & $0.0518$                          \\\cline{2-5}
		                                        & $10^4$  & $10^4$  & $0.0530$                              & $0.0507$                          \\\cline{2-5}
                                                & $10^5$  & $10^5$  & $0.0502$                             & $0.0499$                         \\\hline

\multirow{3}{*}{$n^{1.5}$}                      & $10^3$  & $10^4$  & $0.0456$                              & $0.0526$                          \\\cline{2-5}
		                                        & $10^4$  & $10^4$  & $0.0455$                              & $0.0476$                          \\\cline{2-5}
                                                & $10^5$  & $10^5$  & $0.0484$                             & $0.0496$                         \\\hline

\multirow{3}{*}{$n^{1.8}$}                      & $10^3$  & $10^4$  & $0.0453$                              & $0.0481$                          \\\cline{2-5}
		                                        & $10^4$  & $10^4$  & $0.0485$                              & $0.0489$                          \\\cline{2-5}
                                                & $10^5$  & $10^5$  & $0.0437$                             & $0.0505$                         \\\hline

\multirow{3}{*}{$\frac{n^{2}}{\log(1+n)}$}      & $10^3$  & $10^4$  & $0.0442$                              & $0.0532$                          \\\cline{2-5}
		                                        & $10^4$  & $10^4$  & $0.0426$                              & $0.0514$                          \\\cline{2-5}
                                                & $10^5$  & $10^5$  & $0.0435$                             & $0.0507$                         \\\hline

\end{tabular}}
\end{center}

The table 1 illustrates that  the percentages, $r^o_{N}(b_n)$, are all nearly 0.05 when the truncated sequences $b_n\leq n^{4/3}/10\log(n)$. The simulation results of  $r^o_{N}(b_n)$ are not good for the truncated sequences $b_n=n^{1.5}, n^{1.8}, n^{2}/\log(1+n)$ because the total number of samples and the number of repetitions are not large enough.

\vskip 2cm

\begin{center}
\tabcolsep 0.0000000001in
\setlength{\tabcolsep}{6mm}{
\vskip 0.3cm
{{\bf Table 2:} Simulation  values of  $r_{N}(b_n)$ and $\widetilde{r}_N(n)$ for $\alpha_0=1/2$, $\beta=0.05$  and $y=255.$    \\}
\vskip 0.6cm
\begin{tabular}{|c|c|c|c|c|c|}\hline

$b_n$                                           &$n$      & $N$     & $r_{N}(b_n)$    &$\widetilde{r}_N(n)$    \cr\cline{1-5}
\multirow{3}{*}{$\log(n)$}                      & $10^3$  & $10^4$  & $0.0498$                              & $0.0534$                          \\\cline{2-5}
		                                        & $10^4$  & $10^4$  & $0.0488$                              & $0.0524$                          \\\cline{2-5}
                                                & $10^5$  & $10^5$  & $0.0508$                             & $0.0494$                         \\\hline

\multirow{3}{*}{$n^{0.5}$}                      & $10^3$  & $10^4$  & $0.0527$                              & $0.0472$                          \\\cline{2-5}
		                                        & $10^4$  & $10^4$  & $0.0502$                              & $0.0493$                          \\\cline{2-5}
                                                & $10^5$  & $10^5$  & $0.0504$                             & $0.0497$                         \\\hline

\multirow{3}{*}{$n $}                           & $10^3$  & $10^4$  & $0.0593$                              & $0.0472$                          \\\cline{2-5}
		                                        & $10^4$  & $10^4$  & $0.0557$                              & $0.0513$                          \\\cline{2-5}
                                                & $10^5$  & $10^5$  & $0.0509$                             & $0.0507$                         \\\hline

\multirow{3}{*}{$\frac{n^{4/3}}{10*\log(n)}$}   & $10^3$  & $10^4$  & $0.0578$                              & $0.0500$                          \\\cline{2-5}
		                                        & $10^4$  & $10^4$  & $0.0531$                              & $0.0495$                          \\\cline{2-5}
                                                & $10^5$  & $10^5$  & $0.0519$                             & $0.0494$                         \\\hline

\multirow{3}{*}{$n^{1.5}$}                      & $10^3$  & $10^4$  & $0.0100$                             & $0.0529$                          \\\cline{2-5}
		                                        & $10^4$  & $10^4$  & $0.0827$                              & $0.0476$                          \\\cline{2-5}
                                                & $10^5$  & $10^5$  & $0.0565$                             & $0.0498$                         \\\hline

\multirow{3}{*}{$n^{1.8}$}                      & $10^3$  & $10^4$  & $0.0199$                             & $0.0506$                          \\\cline{2-5}
		                                        & $10^4$  & $10^4$  & $0.0162$                             & $0.0494$                          \\\cline{2-5}
                                                & $10^5$  & $10^5$  & $0.0114$                            & $0.0499$                         \\\hline

\multirow{3}{*}{$\frac{n^{2}}{\log(1+n)}$}      & $10^3$  & $10^4$  & $0.0205$                             & $0.0499$                          \\\cline{2-5}
		                                        & $10^4$  & $10^4$  & $0.0189$                             & $0.0467$                          \\\cline{2-5}
                                                & $10^5$  & $10^5$  & $0.0177$                            & $0.0497$                         \\\hline

\end{tabular}}
\end{center}

Similarly, the table 2 shows that  the percentages, $r_{N}(b_n)$, are all nearly 0.05 when the truncated sequences $b_n\leq n^{4/3}/10\log(n)$. The simulation results of  $r_{N}(b_n)$ are not good for the truncated sequences $b_n=n^{1.5}, n^{1.8}, n^{2}/\log(1+n)$ because the total number of samples and the number of repetitions are not large enough.

\section{ Conclusion}

In order to deals with the hypothesis test for the extremely heavy-tailed distributions with infinite mean or variance, we present two test statistics $T_n$ and $T^o_n$  by using a truncated sample mean according to the known variance and unknown variance. Here we use the sample variance  $\hat{B}^2(b_n) $ to replace the unknown variance. We obtain three  necessary and sufficient conditions in Theorems 2 and 3 for the asymptotic distributions of  $T_n$ and $T^o_n$. 

By $0<h_1\leq D_n(b_n)\leq h_2<\infty$ we can determine a critical  truncated sequence $\{c_n\}$ such that when the truncated sequence is less than or greater than the critical sequence (i.e. $b_n\leq c_n/\widetilde{L}(n)$ or $b_n\geq c_n\widetilde{L}(n)$), the corresponding truncated test statistics $T_n$  converge to normal distribution or $-\infty$ or the combination of stable distributions (   $T^o_n$ converges to normal distribution or 0),  respectively, and when the truncated sequence is equal to the critical sequence, we have $T_{n}\Rightarrow \xi/\sqrt{\eta}$ ( $T^o_n \Rightarrow T_{\alpha, h}$).

Though we can get the characteristic function of the random variable $T_{\alpha, h}$, it is difficult to obtain its density function.  The probability distribution  (NNS-distribution) of $T_{\alpha, h}$  is neither normal nor stable. But there is a close relationship between $T_{\alpha, h}$ and stable random variable $ S_{\alpha, 1}$.

\par
\par

\newpage
\vskip 3cm
\textbf{APPENDIX : Proofs of Lemma 1 and Theorems }
\bigskip
\setcounter{equation}{0}
\bigskip
\renewcommand\theequation{A. \arabic{equation}}

\normalsize

\bigskip

\textbf{Proof of Lemma 1.} By condition (3) we know that condition (I) holds.  Let $0<\alpha_k<1$  and $r\geq 1$. Then
\begin{eqnarray*}
\textbf{E}(X^{r}_k(b_n))&=&\int_{a_k}^{b_n}x^{r}dF_k(x)\\
&=&-b_n^{r}[1-F_k(b_n)]+a_k^{r}[1-F_k(a_k)]+r\int_{a_k}^{b_n}x^{r-1}[1-F_k(x)]dx\\
&=&-b_n^{r-\alpha_k}L(b_n)+a_k^{r-\alpha_k}L(a_k)+r L(\xi_{n, k})\int_{a_k}^{b_n}x^{r-1-\alpha_k}dx\\
&=&b_n^{r-\alpha_k}L(b_n)\Big(\frac{r}{r-\alpha_k}\frac{L(\xi_{n, k})}{L(b_n)}-1\Big)+a_k^{r-\alpha_k}\Big(L(a_k)-\frac{r L(\xi_{n, k})}{r-\alpha_k}\Big)
\end{eqnarray*}
and therefor
\begin{eqnarray}
1=\frac{b_n^{r-\alpha_k}L(b_n)}{\textbf{E}(X^{r}_k(b_n))}\Big(\frac{r}{r-\alpha_k}\frac{L(\xi_{n, k})}{L(b_n)}-1+\frac{a_k^{r-\alpha_k}L(a_k)}{b_n^{r-\alpha_k}L(b_n)}[1-\frac{r L(\xi_{n, k})}{(r-\alpha_k)L(a_k)}]\Big)
\end{eqnarray}
for $1\leq k\leq n$ and large $n$, where the positive number $\xi_{n, k}$ satisfying $a_k< \xi_{n, k} < b_n$, comes from the integral mean value theorem. We can further prove that  $L(\xi_{n, k})/L(b_n) \to 1$ as $n\to \infty$. In fact, let $L(x)\geq L(y)$ for $x\geq y\geq a_k$,  for any small positive number $\varepsilon$ we have
\begin{eqnarray*}
\frac{L(\varepsilon b_n)}{r-\alpha_k}\Big(b_n^{r-\alpha_k}-(\varepsilon b_n)^{r-\alpha_k}\Big)&=& L(\varepsilon b_n)\int_{\varepsilon b_n}^{b_n}x^{r-1-\alpha_k}dx\\
&\leq &\int_{\varepsilon b_n}^{b_n}x^{r-1-\alpha_k}L(x)dx \leq \int_{a_k}^{b_n}x^{r-1-\alpha_k}L(x)dx\\
&= &L(\xi_{n, k})\int_{a_k}^{b_n}x^{r-1-\alpha_k}dx \leq L(b_n) \int_{a_k}^{b_n}x^{r-1-\alpha_k}dx\\
&=& \frac{L(b_n)}{r-\alpha_k}\Big(b_n^{r-\alpha_k}-a_k^{r-\alpha_k}\Big)
\end{eqnarray*}
for $b_n \varepsilon > a_k$, and therefore
\begin{eqnarray*}
\frac{L(\varepsilon b_n)}{L(b_n)}\frac{1-(\varepsilon )^{r-\alpha_k}}{1-(a_k/b_n)^{r-\alpha_k}}\leq \frac{L(\xi_{n, k})}{L(b_n)}\leq 1
\end{eqnarray*}
for large $n$. Thus, by using (4) we have $L(\varepsilon b_n)/L(b_n) \to 1$ and therefore, $L(\xi_{n, k})/L(b_n) \to 1$ as $n\to \infty$. Similarly, we can show that $L(\xi_{n, k})/L(b_n) \to 1$  as $n\to \infty$ when  $L(x)\leq L(y)$ for $x\geq y\geq a_k$. This means that
\begin{eqnarray}
d_{k, n}(r)\triangleq \frac{r L(\xi_{n, k})}{(r-\alpha_k)L(b_n)}-1 \longrightarrow d_k(r)\triangleq \frac{\alpha_k}{r-\alpha_k}
\end{eqnarray}
as $n\to \infty$. Since $b^{r-\alpha_k}_n L(b_n)=b^{r}_n(1-F_k(b_n)) \to \infty$ for $r\geq 1, k\geq 1$ as $n\to \infty$, it follows from (A.1) and (A.2) that  $\textbf{E}(X^{r}_k(b_n))$ satisfies the condition (II) for large $n$.

Let $1\leq \alpha_k<2$.  It is clear that  $b^{r}_n[1-F_k(b_n)]\to \infty $ for $r\geq 2$ as $n\to \infty$. Similarly, it can be proved that (A.1) and (A.2) hold and therefore,  the condition (II)$^{\prime}$ also hold.
\\
\\
\textbf{Proof of Theorem 1.}  (i)  Let $h_n \to \infty$. We first prove that
\begin{eqnarray*}
\frac{T_n}{T^o_n}=[\frac{\hat{B}^2(b_n)}{\sum_{k=1}^n(X_k(b_n)-\mu_0(b_n))^2}\frac{\sum_{k=1}^n(X_k(b_n)-\mu_0(b_n))^2}{\sum_{k=1}^n\textbf{E}(X_k(b_n)-\mu_0(b_n))^2}\frac{\sum_{k=1}^n\textbf{E}(X_k(b_n)-\mu_0(b_n))^2}
{\sum_{k=1}^nVar(X_k(b_n))}]^{-1/2} \to 1
\end{eqnarray*}
in probability as $n\to \infty$. Note that
\begin{eqnarray}
\frac{\sum_{k=1}^n\textbf{E}(X_k(b_n)-\mu_0(b_n))^2}
{\sum_{k=1}^nVar(X_k(b_n))}&=&1+\frac{ \sum_{k=1}^n(\mu_k(b_n)-\mu_0(b_n))^2}{\sum_{k=1}^nVar(X_k(b_n))}\\
\frac{\hat{B}^2(b_n)}{\sum_{k=1}^n(X_k(b_n)-\mu_0(b_n))^2}&=&1-\frac{[\sum_{k=1}^n(X_k(b_n)-\mu_k(b_n))]^2}{n\sum_{k=1}^n(X_k(b_n)-\mu_0(b_n))^2}.
\end{eqnarray}
Let $\varepsilon$ be a small positive number. It follows from condition (I) that
\begin{eqnarray*}
&&\frac{\mu_0(b_n)}{b_n}=(1+o(1))\frac{1}{n}\sum_{k=1}^nd_k(1)[1-F_k(b_n)]\to 0\\
&&\frac{ \sum_{k=1}^n(\mu_k(b_n)-\mu_0(b_n))^2}{\sum_{k=1}^nVar(X_k(b_n))}=(1+o(1))\frac{\sum_{k=1}^nd^2_k(1)[1-F_k(b_n)]^2}{\sum_{k=1}^nd_k(2)[1-F_k(b_n)]}\to 0
\end{eqnarray*}
and therefore, by (A.3)
\begin{eqnarray}
\frac{\sum_{k=1}^n\textbf{E}(X_k(b_n)-\mu_0(b_n))^2}{\sum_{k=1}^nVar(X_k(b_n))} \to 1
\end{eqnarray}
in probability as $n\to \infty$. By  Chebyshev inequality, the condition (II) or (II) $^{\prime}$, and $\sum_{k=1}^n[1-F_k(b_n)] \to \infty$, we can get that
\begin{eqnarray*}
&&\textbf{P}\Big( |\frac{[\sum_{k=1}^n(X_k(b_n)-\mu_k(b_n))]^2}{\sum_{k=1}^n\textbf{E}(X_k(b_n)-\mu_0(b_n))^2}-1|\geq \varepsilon \Big)\leq \frac{\sum_{k=1}^n\textbf{E}(X_k(b_n)-\mu_k(b_n))^4}{ \varepsilon^2[\sum_{k=1}^n\textbf{E}(X_k(b_n)-\mu_0(b_n))^2]^2}\to 0\\
&&\textbf{P}\Big( |\frac{\sum_{k=1}^n(X_k(b_n)-\mu_0(b_n))^2}{\sum_{k=1}^n\textbf{E}(X_k(b_n)-\mu_0(b_n))^2}-1|\geq \varepsilon \Big)\leq  \frac{\sum_{k=1}^n\textbf{E}(X_k(b_n)-\mu_0(b_n))^4}{ \varepsilon^2[\sum_{k=1}^n\textbf{E}(X_k(b_n)-\mu_0(b_n))^2]^2}\to 0
\end{eqnarray*}
as $n\to \infty$. Hence
\begin{eqnarray}
\frac{[\sum_{k=1}^n(X_k(b_n)-\mu_k(b_n))]^2}{\sum_{k=1}^n\textbf{E}(X_k(b_n)-\mu_0(b_n))^2}\to 1,\,\,\,\,\,\,\,\,\, \frac{\sum_{k=1}^n(X_k(b_n)-\mu_0(b_n))^2}{\sum_{k=1}^n\textbf{E}(X_k(b_n)-\mu_0(b_n))^2}\to 1
\end{eqnarray}
in probability as $n\to \infty$ and therefore, by (A.4) and (A.6) we have
\begin{eqnarray}
\frac{\hat{B}^2(b_n)}{\sum_{k=1}^n(X_k(b_n)-\mu_0(b_n))^2}\to 1
\end{eqnarray}
in probability as $n\to \infty$. By (A.5), (A.6) and (A.7), we know that
\begin{eqnarray}
\frac{T_n}{T^o_n}=\frac{[\sum_{k=1}^nVar(X_k(b_n))]^{1/2}}{\hat{B}(b_n)}\to 1
\end{eqnarray}
in probability as $n\to \infty$. Furthermore, since
\begin{eqnarray}
 \frac{\sum_{k=1}^n\textbf{E}|X_k(b_n)-\mu_k(b_n)|^3}{[\sum_{k=1}^nVar(X_k(b_n))]^{3/2}}=O\Big([\sum_{k=1}^n(1-F_k(b_n))]^{-1/2}\Big) \to 0
\end{eqnarray}
as  $n\to \infty$, it follows from  Lyapunov's central limit theorem that $T^o_n \Rightarrow N(0, 1)$  and therefore, by (A.8), $T_n =T^o_nT_n/T^o_n\Rightarrow N(0, 1)$ as  $n\to \infty$.

(ii). Let $h_n \to h$. Let $X_{nk}=X_k(b_n)/b_n$, $\mu_{nk}=\mu_k(b_n)/b_n$ and $\mu_{0n}=\mu_0(b_n)/b_n$. Let $F_n(x, y)$ and $C_n(t_1, t_2)$ denote respectively  the distribution function and the characteristic function of $(\xi_n, \eta_n)$ for $n\geq 1$. Since the observation sequence $\{X_k\}$ is mutually independent, it follows from conditions (I), (II) or (II)$^{\prime}$ that $\mu_{nk}=\mu_k(b_n)/b_n \to 0$, $\mu_{0n}=\mu_0(b_n)/b_n \to 0$ as $n\to \infty$ and for large $n$,
\begin{eqnarray*}
&&\sum_{k=1}^n\textbf{E}([t_1(X_{nk}-\mu_{nk})+t_2(X_{nk}-\mu_{0n})^2])=(1+o(1))h_n(2)t_2\\
&&\sum_{k=1}^n\textbf{E}([t_1(X_{nk}-\mu_{nk})+t_2(X_{nk}-\mu_{0n})^2]^m)=(1+o(1))\sum_{j=0}^mh_n(2m-j)C^j_m t^j_1t^{m-j}_2
\end{eqnarray*}
for $m\geq 2$ and
\begin{eqnarray*}
\sum_{m=N+1}^{\infty}\sum_{k=1}^n\textbf{E}([t_1(X_{nk}+\mu_{nk})+t_2(X_{nk}+\mu_{0n})^2]^m)/m!=o(1)
\end{eqnarray*}
for large $N$, where $N$ is large such that  $(|t_1|+|t_2|)/N$ is small.  Since $h_n \to h$, it follows from (6), the conditions (II) and (II)$^{\prime}$ that
\begin{eqnarray*}
d_{n}(m)\to \alpha/(m-\alpha)\,\,\,\,\text{ and}\,\,\,\, h_{n}(m) \to h(m)=h\alpha/(m-\alpha)
\end{eqnarray*}
as $n\to \infty$, where $0<\alpha <1$ or $1\leq \alpha <2$ respectively  under the condition (II) or the condition (II)$^{\prime}$. Hence, the limiting function of characteristic  functions   of $(\xi_{n}, \eta_{n})$ can be written as
\begin{eqnarray*}
C(t_1, t_2)&=&\lim_{n\to\infty}C_{n}(t_1, t_2)=\lim_{n\to\infty}\textbf{E}\Big(e^{i(t_1\xi_{n}+t_2\eta_{n})}\Big)\\
&=&\lim_{n\to \infty}\prod_{k=1}^{n}\textbf{E}\Big(\exp\{i[t_1(X_{nk}-\mu_{nk})+t_2(X_{nk}-\mu_{0n})^2]\}\Big)\\
&=&\lim_{n\to \infty}\exp\{(1+o(1))\Big(ih_{n}(2)t_2+\sum_{m=2}^{\infty}\frac{i^m}{m!}\sum_{j=0}^mh_{n}(2m-j)C^j_m t^j_1t^{m-j}_2\Big)\}\\
&=&\exp\{ih(2)t_2+\sum_{m=2}^{\infty}i^m\sum_{j=0}^m\frac{h(2m-j)t^j_1t^{m-j}_2}{j!(m-j)!}\}.
\end{eqnarray*}
Since $C(t_1, t_2)$ is continuous at $(t_1, \, t_2)=(0, \,0)$, it follows that there is a two-dimensional random vector $(\xi, \, \eta)$ with the distribution function $F(x, y)$ corresponding to the characteristic function $C(t_1, t_2)$ such that $(\xi_{n}, \,  \eta_{n}) \Rightarrow (\xi,\, \eta)$ as $n\to \infty$.

Next we prove that $\xi_{n}/\sqrt{\eta_{n}} \Rightarrow \xi/\sqrt{\eta}$  as $n\to \infty$.
For any small $\varepsilon>0$, we can take two  large $M$ and $N_1$ such that
\begin{eqnarray*}
&&|F_{n}(x, +\infty)-F_{n}(x, M)-[F(x, +\infty)-F(x, M)]|\\
&&+|F_{n}(+\infty, y)-F_{n}(M, y)-[F(+\infty, y)-F(M, y)]|+|F_{n}(-M, y)-F(-M, y)]|<\varepsilon
\end{eqnarray*}
for $n\geq N_1$, where $(x, M)$, $(-M, y)$  and $ (M,  y)$ are all continuous points of $F(.,.)$. By the properties of distribution function, we need only consider $|x|\leq M$ and $0\leq y\leq M$.

For any fixed real number $z$, let $I(x, y)=I(x/\sqrt{y}\leq z)$  be a indicator function and $g_m(x, y)$ are bounded continuous functions satisfying $ g_m(x, y)\nearrow I(x, y)$ as $m \to \infty$ for $|x|\leq M$ and $0\leq y\leq M$. It follows that
\begin{eqnarray*}
&&|\textbf{P}\Big( \frac{\xi_{n}}{\sqrt{\eta_{n}}} \leq z\Big)-\textbf{P}\Big( \frac{\xi}{\sqrt{\eta}} \leq z\Big)|\\
&&=|\int_{}^{}\int_{}^{}I(x, y)dF_{n}(x, y)-\int_{}^{}\int_{}^{}I(x, y)dF(x, y)|\\
&&\leq \int_{}^{}\int_{}^{}(I(x, y)-g_m(x, y))dF_{n}(x, y)+\int_{}^{}\int_{}^{}(I(x, y)-g_m(x, y))dF(x, y)\\
&&+|\int_{}^{}\int_{}^{}g_m(x, y)dF_{n}(x, y)-\int_{}^{}\int_{}^{}g_m(x, y)dF(x, y)| \, \to \,0
\end{eqnarray*}
as $n\to \infty$, where the second term of the last inequality comes from the equivalent property of convergence in distribution. This means that $\xi_{n}/\sqrt{\eta_{n}} \Rightarrow \xi/\sqrt{\eta}$  as $n\to \infty$.

Let $t_2=0$ or $t_1=0$ in $C(t_1, t_2)$, we can get the characteristic functions of  $\xi$ and $\eta$,   $C_{\xi}(t)=\exp\{iu_1(t)+v_1(t)\}$ and $C_{\eta}(t)=\exp\{iu_2(t)+v_2(t)\}$, respectively, where
\begin{eqnarray*}
&&u_1(t)=\alpha h\sum_{k=1}^{\infty}\frac{(-1)^k}{(2k+1-\alpha)}\frac{t^{2k+1}}{(2k+1)!},\,\,\,\,\,\,\,v_1(t)=\alpha h\sum_{k=1}^{\infty}\frac{(-1)^k}{(2k-\alpha)}\frac{t^{2k}}{(2k)!}\\
&&u_2(t)=\frac{\alpha h}{2}\sum_{k=0}^{\infty}\frac{(-1)^k}{(2k+1-\alpha/2)}\frac{t^{2k+1}}{(2k+1)!},\,\,\,\,\,\,\,v_2(t)=\frac{\alpha h}{2}\sum_{k=1}^{\infty}\frac{(-1)^k}{(2k-\alpha/2)}\frac{t^{2k}}{(2k)!}
\end{eqnarray*}
By the expansion of power series of $sinx$ and $cosx$, we can get that
\begin{eqnarray*}
&&u_1(t)=\alpha h |t|^{\alpha}\, sgn(t)\int_{0}^{|t|}\frac{sinx-x}{x^{1+\alpha}}dx,\,\,\,\,\,\,\,\, v_1(t)=-\alpha h |t|^{\alpha}\int_{0}^{|t|}\frac{1-cosx}{x^{1+\alpha}}dx\\
&&u_2(t)=\frac{\alpha h}{2} |t|^{\frac{\alpha}{2}}\, sgn(t)\int_{0}^{|t|}\frac{sinx}{x^{1+\frac{\alpha}{2}}}dx,\,\,\,\,\,\,\,\,v_2(t)=-\frac{\alpha h}{2} |t|^{\frac{\alpha}{2}}\int_{0}^{|t|}\frac{1-cosx}{x^{1+\frac{\alpha}{2}}}dx
\end{eqnarray*}
for $t\in (-\infty,\, \infty)$. Hence, (7) and (8) hold for $0<\alpha <2$.

(iii) Let $h_n \to 0$ and $0<\alpha<1$. Replacing $h$ in (7) and (8) with $h_n$, we can get the characteristic functions of  $\xi_n$ and $\eta_n$,   $C_{\xi_n}(t)=\exp\{iu_{1n}(t)+v_{1n}(t)\}$ and $C_{\eta_n}(t)=\exp\{iu_{2n}(t)+v_{2n}(t)\}$, respectively, where
\begin{eqnarray*}
&&u_{1n}(t)=(1+o(1))\alpha h_n |t|^{\alpha}\, sgn(t)\int_{0}^{|t|}\frac{sinx-x}{x^{1+\alpha}}dx,\,\,\,\,   v_{1n}(t)=-(1+o(1))\alpha h_n |t|^{\alpha}\int_{0}^{|t|}\frac{1-cosx}{x^{1+\alpha}}dx\\
&&u_{n2}(t)=(1+o(1))\frac{\alpha h_n}{2} |t|^{\frac{\alpha}{2}}\, sgn(t)\int_{0}^{|t|}\frac{sinx}{x^{1+\frac{\alpha}{2}}}dx,\,\,\,\,v_{2n}(t)=-(1+o(1))\frac{\alpha h_n}{2} |t|^{\frac{\alpha}{2}}\int_{0}^{|t|}\frac{1-cosx}{x^{1+\frac{\alpha}{2}}}dx
\end{eqnarray*}
for large $n$. Note that $C_{\xi_nh^{-\frac{1}{\alpha}}_n}(t)=C_{\xi_n}(th^{-\frac{1}{\alpha}}_n)$, $C_{\eta_n h^{-\frac{2}{\alpha}}_n}(t)=C_{\eta_n}(th^{-\frac{2}{\alpha}}_n)$,
\begin{eqnarray}
u_{1n}(th^{-\frac{1}{\alpha}}_n)&=&(1+o(1))\alpha |t|^{\alpha}\, sgn(t)\int_{0}^{|t|h^{-\frac{1}{\alpha}}_n}\frac{sinx-x}{x^{1+\alpha}}dx\\
&=&(1+o(1))\alpha |t|^{\alpha}\, sgn(t)\Big(\int_{0}^{1}\frac{sinx-x}{x^{1+\alpha}}dx+\int_{1}^{|t|h^{-\frac{1}{\alpha}}_n}\frac{sinx}{x^{1+\alpha}}dx-\frac{((|t|h^{-\frac{1}{\alpha}}_n)^{1-\alpha}-1)}{1-\alpha}\Big)\nonumber \\
&\longrightarrow & -\infty\nonumber
\end{eqnarray}
and
\begin{eqnarray}
&&v_{1n}(th^{-\frac{1}{\alpha}}_n) \to -\alpha |t|^{\alpha}\int_{0}^{\infty}\frac{1-cosx}{x^{1+\alpha}}dx,\\
&&v_{2n}(th^{-\frac{1}{\alpha}}_n)\to -\frac{\alpha }{2} |t|^{\frac{\alpha}{2}}\int_{0}^{\infty}\frac{1-cosx}{x^{1+\frac{\alpha}{2}}}dx,\\
&&u_{n2}(th^{-\frac{1}{\alpha}}_n) \to \frac{\alpha }{2} |t|^{\frac{\alpha}{2}}\, sgn(t)\int_{0}^{\infty}\frac{sinx}{x^{1+\frac{\alpha}{2}}}dx
\end{eqnarray}
as $n\to \infty$. Hence, $\xi_nh^{-\frac{1}{\alpha}}_n \to -\infty$ in probability and $\eta_n h^{-\frac{2}{\alpha}}_n$ converges to a finite random variable. Since $\hat{B}^2(b_n)/b^2_n=\eta_n-\xi^2_n/n$, it follows that
\begin{eqnarray*}
T_{n}=\frac{\xi_{n}}{\sqrt{\eta_{n}}}\Big(1-\frac{\xi^2_{n}}{n\eta_{n}}\Big)^{-1/2}=\frac{\xi_{n}h^{-\frac{1}{\alpha}}_n}{\sqrt{\eta_{n}h^{-\frac{2}{\alpha}}_n}}\Big(1-\frac{\xi^2_{n}h^{-\frac{2}{\alpha}}_n}{n\eta_{n}h^{-\frac{2}{\alpha}}_n}\Big)^{-1/2} \to -\infty
\end{eqnarray*}
as $n\to \infty$ in probability.

Let  $h_n \to 0$ and $\alpha=1$. Similarly, we have
\begin{eqnarray}
u_{1n}(th^{-1}_n)&=&(1+o(1)) |t|\, sgn(t)\int_{0}^{|t|h^{-1}_n}\frac{sinx-x}{x^{2}}dx\\
&=&(1+o(1))|t|\, sgn(t)\Big(\int_{0}^{1}\frac{sinx-x}{x^{2}}dx+\int_{1}^{|t|h^{-1}_n}\frac{sinx}{x^{2}}dx-\log(|t|h^{-1}_n)\Big)\nonumber\\
&\longrightarrow & -\infty\nonumber
\end{eqnarray}
as $n\to -\infty$ and therefore, $\xi_nh^{-1}_n \to -\infty$ in probability and $\eta_n h^{-2}_n$ converges to a finite random variable. Thus, $ T_{n} \to -\infty$ as $n\to \infty$.

(iv). Let $h_n \to 0$ and $1<\alpha<2$. Similar to (A.10)-(A.13), we can get that the limiting characteristic  functions, $C_{\widetilde{\xi}}(t)$ and $C_{\widetilde{\eta}}(t)$ of $(\xi_{n}h_n^{-1/\alpha},\,\,\eta_{n}h_n^{-2/\alpha})$ in the following
\begin{eqnarray*}
C_{\widetilde{\xi}}(t)&=&\exp\{\alpha |t|^{\alpha}\Big(i\, sgn(t)\int_{0}^{\infty}\frac{sinx-x}{x^{1+\alpha}}dx-\int_{0}^{\infty}\frac{1-cosx}{x^{1+\alpha}}dx\Big)\}\\
&=&\exp\{\alpha |t|^{\alpha}\Big(\int_{0}^{\infty}\frac{e^{isgn(t)x}-1-isgn(t)x}{x^{1+\alpha}}dx\Big)\}\\
C_{\widetilde{\eta}}(t)&=&\exp\{\frac{\alpha}{2} |t|^{\frac{\alpha}{2}}\Big(i\, sgn(t)\int_{0}^{\infty}\frac{sinx}{x^{1+\frac{\alpha}{2}}}dx-\int_{0}^{\infty}\frac{1-cosx}{x^{1+\frac{\alpha}{2}}}dx\Big)\}\\
&=&\exp\{\frac{\alpha}{2} |t|^{\frac{\alpha}{2}}\Big(\int_{0}^{\infty}\frac{e^{i sgn(t)x}-1}{x^{1+\frac{\alpha}{2}}}dx\Big)\}
\end{eqnarray*}
respectively. Furthermore, by the characteristic  function of $\Gamma$ distribution, we have (see [6])
\begin{eqnarray*}
C_{\widetilde{\xi}}(t)&=&\exp\{-c_1 |t|^{\alpha}(1-isgn(t)tan(\frac{\alpha \pi}{2}))\},\\
C_{\widetilde{\eta}}(t)&=&\exp\{-c_2|t|^{\frac{\alpha}{2}}(1-isgn(t)tan(\frac{\alpha \pi}{4}))\},
\end{eqnarray*}
where $c_1=\Gamma(2-\alpha)|cos(\alpha \pi/2)|/(\alpha-1)$ and $c_2=2\Gamma(1-\alpha/2)cos(\alpha \pi/4)$. It is clare that $C_{\widetilde{\xi}}(t)$ and $C_{\widetilde{\eta}}(t)$ are the characteristic  functions of stable random variables $S_{\alpha, -1}$ and $S_{\alpha/2, 1}$ respectively. It completes the proof.\\
\\

\textbf{Proof of Theorem 2.} We only need to prove the necessity. (i) Let $T_n \Rightarrow N(0, 1)$. Assume that $h_n$ dose not converge to $\infty$, that is, there is s subsequence $\{n_l\}$ such that  $h_{n_l}\to h$, where $0\leq h<\infty$. By (ii)-(iv) of Theorem 1, we know that $ T_{n_l} \Rightarrow \xi/\sqrt{\eta}$ or $T_{n_l} \to -\infty$ in probability or $T_{n_l} \Rightarrow S_{\alpha}/\sqrt{S_{\alpha/2}}$. Note that $\xi/\sqrt{\eta}$ and $S_{\alpha, -1}/\sqrt{S_{\alpha/2, 1}}$ are not standard normal since the two distribution functions of $\xi$ and $S_{\alpha, -1}$ are not symmetric. This contradiction means that $h_n\to \infty$.

(ii) Let $ T_{n} \Rightarrow \xi/\sqrt{\eta}$. Assume that there is s subsequence $\{n_l\}$ such that $h_{n_l} \to \infty$ or $h_{n_l} \to 0$. By (i), (iii) and (iv) of Theorem 1, we know that $T_{n_l} \Rightarrow N(0, 1)$ or $T_{n_l} \to -\infty$ in probability or $T_{n_l} \Rightarrow S_{\alpha, -1}/\sqrt{S_{\alpha/2, 1}}$. This contradiction means that $h_{n} \to h$, where $0<h<\infty$.

(iii) Let $T_{n} \to -\infty$.  Assume that there is s subsequence $\{n_l\}$ such that $h_{n_l} \to \infty$ or $h_{n_l} \to h>0$. By (i) and (ii) of Theorem 1, we see that $T_{n_l} \Rightarrow N(0, 1)$ or $ T_{n_l} \Rightarrow \xi/\sqrt{\eta}$. This contradiction means that $h_{n} \to 0$ and therefore, $0<\alpha\leq 1$ since $T_{n} \to -\infty$.

 (iv) Let $T_{n} \Rightarrow S_{\alpha, -1}/\sqrt{S_{\alpha/2, 1}}$. Similar to the proof of (iii), we have $h_{n} \to 0$ and $1<\alpha <2$. It completes the proof.
\\

\textbf{Proof of Theorem 3.} Like the proof of Theorem 1, we first prove the sufficiency of Theorem 3. Let $h_n\to \infty$. By  (A.9) and Lyapunov's central limit theorem, we know that $T^o_n \Rightarrow N(0, 1)$.

Let $h_n\to h, 0<h<\infty$. Note that $\sqrt{Var(b_n)}=(1+o(1))\sqrt{h}\sigma b_n$ for large $n$. Replacing $t$ in (7) with $t/\sqrt{h}\sigma$, we can get the limiting characteristic  function $C_{\alpha, h}(t)$ of the characteristic  functions $C_{ T^o_n}(t)$ for $T^o_n, n\geq 1$, that is,
\begin{eqnarray*}
C_{\alpha, h}(t)=\exp\{h\alpha \sigma^{-\alpha} |t|^{\alpha}\Big(i\, sgn(t)\int_{0}^{|t|/\sqrt{h}\sigma}\frac{sinx-x}{x^{1+\alpha}}dx-\int_{0}^{|t|/\sqrt{h}\sigma}\frac{1-cosx}{x^{1+\alpha}}dx\Big)\}
\end{eqnarray*}
for $0<\alpha <2$. This is (11).

Let  $h_n\to 0$ and $0<\alpha<1$. Since $\sqrt{Var(b_n)}=(1+o(1))\sigma\sqrt{h_n} \,b_n$, it follows from (5), (6) and (A.10) that
\begin{eqnarray*}
T^o_n&=& \frac{\sum_{k=1}^n(X_k(b_n)-\mu_k(b_n)  )}{\sqrt{Var(b_n)}}=(1+o(1))\sigma^{-1}h_n^{\frac{1}{\alpha}-\frac{1}{2}}\xi_n h_n^{-\frac{1}{\alpha}}\\
&=&O\Big(h_n^{\frac{1}{\alpha}-\frac{1}{2}}h_n^{-\frac{1-\alpha}{\alpha}}\Big)\to 0
\end{eqnarray*}
as $n\to \infty$, since $(\xi_n h_n^{-\frac{1}{\alpha}})h_n^{\frac{1-\alpha}{\alpha}} $  converges to $-\alpha/(1-\alpha)$ in probability.\\
Let  $h_n\to 0$ and $\alpha=1$. Similarly, by (A.14) we have $T^o_n=O(h_n^{\frac{1}{2}}\log h_n) \to 0$ as $n\to \infty$.\\
Let  $h_n\to 0$ and $1<\alpha <2$. From (iv) of Theorem 1 it follows that
\begin{eqnarray*}
T^o_n=(1+o(1))\sigma^{-1}h_n^{\frac{1}{\alpha}-\frac{1}{2}}(\xi_n h_n^{-\frac{1}{\alpha}}) \to 0
\end{eqnarray*}
as $n\to \infty$.

By the similar method of proving Theorem 2, we can prove the necessity of Theorem 3. It completes the proof.
\\
\\
\textbf{Proof of Theorem 4.} It needs only to prove that $\sum_{k=1}^n[1-F_k(c_{nk}L^{-1}_0(c_{nk}))] \to \infty$ and $\sum_{k=1}^n[1-F_k(c_{nk}L_0(c_{nk}))] \to 0$ as $n\to \infty$. In fact, since $1-F_k(c_{nk})=h/n$, $c_{nk}=n^{1/\alpha_k}[L(c_{nk})/h]^{1/\alpha_k}$ and therefore
\begin{eqnarray*}
1-F_k(c_{nk}L^{-1}_0(c_{nk}))=(c_{nk}L^{-1}_0(c_{nk}))^{-\alpha_k}L(c_{nk}L^{-1}_0(c_{nk}))=\frac{h}{n}\frac{L_0^{\alpha_k}(c_{nk})L(c_{nk}L^{-1}_0(c_{nk}))}{L(c_{nk})},
\end{eqnarray*}
it follows that $\sum_{k=1}^n[1-F_k(c_{nk}L^{-1}_0(c_{nk}))] \to \infty$ as $n\to \infty$. Similarly, we can prove that $\sum_{k=1}^n[1-F_k(c_{nk}L_0(c_{nk}))] \to 0$ as $n\to \infty$.
\\
\\

\textbf{Proof of Theorem 5.} Note that $h_n\equiv 1$, $\sqrt{Var(c_n)}=(1+o(1))\sigma \,c_n$ and $\sum_{k=1}^n\textbf{E}(X_k(c_n))=(1+o(1))c_n\alpha/(1-\alpha)$, where $\sigma=\sqrt{\alpha}/\sqrt{2-\alpha}$. It follows that
\begin{eqnarray*}
U_n= (1+o(1))\Big(\sigma T^o_n+\frac{\alpha}{1-\alpha }\Big)
\end{eqnarray*}
for $0<\alpha<1$ and
\begin{eqnarray*}
U_n = (1+o(1))\sigma T^o_n
\end{eqnarray*}
for $1\leq \alpha<2$. Thus, from (ii) of Theorem 3 it follows that the results of Theorem 5 hold.
\\
\\

\end{document}